\documentclass[a4paper]{amsart}
\usepackage{amssymb}
\newtheorem{theorem}{Theorem}
\newtheorem{corollary}{Corollary}

\theoremstyle{definition}
\newtheorem*{definition}{Definition}

\begin{document}
\title[Groups in Which Unconditionally Closed Sets are
Algebraic]{A Class of Groups in Which All Unconditionally Closed Sets are
Algebraic}
\author{Ol'ga V. Sipacheva}
\begin{abstract}
It is proved that, in certain subgroups of  direct products of
countable groups, the property of being an unconditionally closed set
coincides with that of being an algebraic set.
In particular, these properties coincide in all Abelian groups.
\end{abstract}
\thanks{This work was financially supported by
the Russian Foundation for Basic Research (project no.~06-01-00764).}
\address
{Department of General Topology and Geometry\\
Mechanics and Mathematics Faculty\\
Moscow State University\\
Leninskie Gory\\
Moscow, 119992 Russia}


\email
{o-sipa@yandex.ru}

\maketitle

Markov~\cite{Markov1945}
called a subset $A$ of a group $G$
\emph{unconditionally closed} in $G$ if it is closed in
any Hausdorff group topology on $G$.  Clearly, all solution sets of
equations in $G$, as well as their finite unions and arbitrary
intersections, are unconditionally closed.

\begin{definition}[Markov~\cite{Markov1945}]
A subset $A$ of a group $G$ with identity element $1$ is said to be
\emph{elementary algebraic} in $G$ if there exists a word $w= w(x)$ in the
alphabet $G\cup \{x^{\pm1}\}$ ($x$ is a variable) such that
$A =\{x\in G: w(x) = 1\}$.
Finite unions of elementary algebraic sets are \emph{additively
algebraic} sets. An arbitrary
intersection of additively algebraic sets is called \emph{algebraic}.
Thus, the algebraic sets in $G$ are the solution sets of arbitrary
conjunctions of finite disjunctions of equations.
\end{definition}

In his 1945 paper~\cite{Markov1945}, Markov showed that any
algebraic set is unconditionally closed and
posed the problem of whether the converse is true.
In~\cite{Markov1946} (see also~\cite{Markov1944}), he solved this
problem for countable groups by proving that any unconditionally
closed set in a countable group is algebraic.
In this paper, we show that the answer is also positive
for some subgroups of direct products of countable groups (in particular,
for all Abelian groups); the proof is based on Markov's ideas.
In the general case, the answer is negative
(at least, under the continuum hypothesis)~\cite{Sipa}.

\begin{definition}[Markov~\cite{Markov1946}]
Let $m$ be a positive integer. By a \emph{multiplicative function} of
$m$ arguments we mean an arbitrary word on the alphabet $\{(1,\pm1),
\dots, (m, \pm1)\}$. The \emph{length} of the function is equal to
the length of this word. Suppose that $G$ is a group and
$x_1, \dots, x_m\in G$. The \emph{value} of a multiplicative function
$\Phi=(j_1, \varepsilon_1)\dots (j_n, \varepsilon_n)$ (of $m$
arguments) at $x_1, \dots, x_n$ in $G$ is defined as
$$
\Phi_G(x_1, \dots, x_n) =\prod_{i=1}^n x_{j_i}^{\varepsilon_i}.
$$
\end{definition}

\begin{definition}
For a subset $A$ of a group $G$, $\mathop{\widetilde A}\nolimits^G$
(or simply $\widetilde A$) denotes the \emph{algebraic closure} of
$A$ in $G$, i.e., the least algebraic set in $G$ containing $A$.
\end{definition}

In what follows, we consider algebraic closures only in $H$.

Following Dikranjan and Shakhmatov, we say that a subgroup $H'$ of a
group $H$ is \emph{supernormal} in $H$ if, for every $x\in H$, there exists
a $y \in H'$ such that $x^{-1} hx = y^{-1} hy$ for all $h\in H$.

By a \emph{countable subproduct} of a direct product
$\prod_{\alpha\in I} G_\alpha$ of groups we mean any product of the form
$\prod_{\alpha\in I} G'_\alpha$, where $G'_\alpha = G_\alpha$ for countably
many indices $\alpha$ and $G'_\alpha= \{1_\alpha\}$ for the remaining
$\alpha$ ($\{1_\alpha\}$ denotes the identity element of the group $G_\alpha$).

\begin{theorem}
If $H$ is a subgroup of a direct product $G$ of countable groups such that
$H\cap G'$
is supernormal in $H$ for any countable subproduct $G'$ of $G$,
then $A\subset H$ is unconditionally closed in $H$ if and only
if it is algebraic in $H$.
\end{theorem}

\begin{proof}
Let $G$ be the direct product of countable groups
$G_\alpha$, where $\alpha\in I$, and let $H$ be a subgroup of $G$ with the
property specified in the statement of the theorem.
By $1$ and
$1_\alpha$ we denote the identity elements of $G$ and $G_\alpha$,
respectively.

Suppose that $A\subset H$, $1\in \widetilde A$, and
$1\notin A$ (in particular, $A$ is not algebraic in $H$).
Let us show that $A$ is not unconditionally closed in $H$.
We set $a_0= 1$ and take an arbitrary element $a_1\in H\setminus \{1\}$.
This element has only finitely many nonzero coordinates; we denote
their indices by $\beta_1, \dots, \beta_{k_1}$.
Note that the set of all multiplicative functions of length
$<6$ with $2$ arguments is finite. Let $\mathfrak M_1$ be the
(finite) subset of this set consisting of all functions $\Phi$
for which
$$ \Phi(a_1, 1) \ne 1.
$$
We set
$$
A_{\Phi}=\{x\in H:\Phi(a_1, x)=1\} \qquad \text{for} \quad \Phi \in
\mathfrak M_1
$$
and
$$
B_1=\bigcup_{\Phi \in \mathfrak M_1}A_\Phi.
$$
The set $B_1$ is additively algebraic.
We have $1\notin B_1$. On the other hand, $1\in \widetilde A$. Since
$\widetilde A$ is the algebraic closure of $A$, we have
$A\setminus B_1\ne \varnothing$. Take $x_1\in A\setminus B_1$.
Note that $x_1\ne a_0, a_1$. Indeed, $x_1\ne a_0= 1$
because $x_1\in A\ni 1$, and $x\ne a_1$ because
$\Phi=(1,1)(1,-1)$ is a multiplicative function of length~2 with
two arguments for which $\Phi(a_1, a_1) = a_1a_1^{-1} = 1$ but
$\Phi(a_1, 1) = a_1\ne 1$.

Let $\{\alpha_1, \dots, \alpha_{n_1}\}$ be a finite
subset of the index set $I$ such that $\{\beta_1, \dots,
\beta_{k_1}\}\subset \{\alpha_1, \dots, \alpha_{n_1}\}$ and
the indices of all nonidentity coordinates of the element $x_1\in
G=\prod_{\alpha\in I} G_\alpha$ are contained in $\{\alpha_1, \dots,
\alpha_{n_1}\}$ (here and in what follows, $\prod$ denotes
direct product).  Let us number all elements $a$ of
the at most countable set $H\cap
\prod_{\alpha\in I} G'_\alpha\setminus \{a_0, a_1\}$, where
$G'_\alpha=G_\alpha$ for $\alpha\in \{\alpha_1, \dots,
\alpha_{n_1}\}$ and $G'_\alpha=\{1_\alpha\}$ for
$\alpha\notin \{\alpha_1, \dots, \alpha_{n_1}\}$,
by positive integers $\ge 2$ so that infinitely many
numbers remain unoccupied but the number $2$ is occupied, i.e.,
$a_2$ is defined (e.g., we can  put $a_2=x_1$ and use
only even numbers for numbering; in this case, we have
$H\cap \prod_{\alpha\in I} G'_\alpha=\{a_0, a_1, a_2, a_4, a_6,
\dots\}$. The set
$H\cap \prod_{\alpha\in I} G'_\alpha$ may be finite, but it necessarily
contains $a_0$, $a_1$, and $a_2$).

Take an element $a_2$. Arguing as above, we choose
$x_2\in A$ so that if
$\Phi$ is a multiplicative function of length $<9$ with $4$
arguments and $\Phi(a_1, a_2, x_1, x_2) = 1$, then $\Phi(a_1, a_2,
x_1, 1) = 1$, and show that $x_2\ne a_0, a_1, a_2$.
Let $n_2 > n_1$, and let $\{\alpha_1, \dots,
\alpha_{n_2}\}\subset I$ be a set of indices which contains the indices
of all nonidentity coordinates of $x_2$.
Let us number all unnumbered
elements $a$ of the at most countable
set $H\cap \prod_{\alpha\in I} G''_\alpha$, where
$G''_\alpha=G_\alpha$ for $\alpha\in \{\alpha_1, \dots,
\alpha_{n_2}\}$ and $G''_\alpha=\{1_\alpha\}$ for
$\alpha\notin \{\alpha_1, \dots, \alpha_{n_2}\}$,
by  positive integers that have not been used at the
preceding step so that infinitely many numbers remain
unoccupied but the number $3$ is occupied (e.g., we can put
$a_3=x_2$ and use
only odd numbers divisible by $3$ for numbering).

At the $j$th step, the situation is as follows.
Positive integers $n_1 < \dots <n_{j-1}$ and indices $\alpha_1, \dots,
\alpha_{n_{j-1}}\in I$  are defined;
elements $x_1, \dots, x_{j-1}\in A$ and
$a_1, \dots, a_j\in H\cap \prod_{\alpha\in I} G'''_\alpha$, where
$G'''_\alpha=G_\alpha$ for $\alpha\in \{\alpha_1, \dots,
\alpha_{n_{j-1}}\}$ and $G''_\alpha=\{1_\alpha\}$ for
$\alpha\notin \{\alpha_1, \dots, \alpha_{n_{j-1}}\}$, are defined; all
elements $a$ of the at most countable
set $H\cap \prod_{\alpha\in I} G'''_\alpha$
are numbered by positive integers, and infinitely many
positive integers remain unoccupied (however, the numbers
$1, \dots, j$ are occupied by the elements $a_1, \dots, a_j$).
Moreover, for any multiplicative function $\Phi$ of length $<3j$
with $2(j-1)$ arguments satisfying the condition $\Phi(a_1, \dots,
a_{j-1}, x_1, \dots, x_{j-1})= 1$, we have $\Phi(a_1, \dots, a_{j-1},
x_1, \dots, x_{j-2}, 1)= 1$.  Let us define $x_j$.

The set of all multiplicative functions of length $<3(j+1)$ with $2j$
arguments is finite. Let $\mathfrak M_j$ be its subset
consisting of all functions $\Phi$
for which
$$
\Phi(a_1,\dots, a_j, x_1, \dots, x_{j-1}, 1) \ne 1.
$$
We set
$$
A_{\Phi}=\{x:\Phi(a_1,\dots, a_j, x_1, \dots, x_{j-1}, x)=1\}
\qquad \text{for} \quad \Phi \in \mathfrak M_j
$$
and
$$
B_j=\bigcup_{\Phi \in \mathfrak M_j}A_\Phi.
$$
The set $B_j$ is additively algebraic, and
$1\notin B_j$. On the other hand, $1\in \widetilde A$. Since
$\widetilde A$ is the algebraic closure of $A$, we have
$A\setminus B_j\ne \varnothing$. Choose any $x_j\in A\setminus B_j$.
Clearly, $x_j\ne
a_0=1, a_1, \dots, a_j$, because if $x_j=a_i$, then $\Phi=(i,
1)(2j, -1)$ is a  multiplicative function of length $< 3(j+1)$ with $2j$
arguments for which $\Phi (a_1, \dots, a_j, x_1, \dots, x_j)=a_ix_j^{-1}=1$
but $\Phi(a_1,\dots, a_j, x_1, \dots, x_{j-1}, 1)\ne1$.
Take sufficiently large number $n_j>n_{j-1}$ and
set $\{\alpha_{n_{j-1}+1}, \dots, \alpha_{n_j}\}\subset I$ for
which all coordinates of the element $x_j$ with indices not belonging to
$\{\alpha_1, \dots, \alpha_{n_j}\}$ are the identity elements of the
corresponding factors.  Let us number all elements $a$ of
the at most countable set $H\cap
\prod_{\alpha\in I} G''''_\alpha$ (where $G''''_\alpha=G_\alpha$ for
$\alpha\in \{\alpha_1, \dots, \alpha_{n_j}\}$ and $G''''_\alpha=\{1_\alpha\}$
for $\alpha\notin \{\alpha_1, \dots, \alpha_{n_j}\}$) which have not
been numbered at the preceding steps by positive integers not occupied
at the preceding steps so that infinitely many numbers remain
unoccupied. If the element
$a_{j+1}$ has not yet been defined (i.e., the number $j+1$ is unoccupied),
then we set $a_{j+1}=x_j$.

As a result, we obtain a countable set $I^*=\{\alpha_i:i\in
\mathbb N\}\subset I$; all elements of the set $G^*=H\cap
\prod_{\alpha\in I} G^*_\alpha$, where
$G^*_\alpha=G_\alpha$ for $\alpha\in I^*$ and $G^*_\alpha=\{1_\alpha\}$
for $\alpha\notin I^*$, are numbered by nonnegative
integers: $G^*=\{a_0, a_1, a_2, \dots\}$, where and $a_0=1$. We
also have a sequence $\{x_1, x_2, \dots\}\subset A$
such that if $\Phi$ is a multiplicative function of length $<3(j+1)$
with $2j$ arguments and $\Phi(a_1,
\dots, a_j, x_1, \dots, x_j)= 1$, then
$\Phi(a_1, \dots, a_j, x_1, \dots, x_{j-1}, 1)= 1$. We set
$A^*=A\cap G^*$.

Note that $G^*$ is a countable supernormal subgroup of $H$ and
$A^*\subset G^*$. Thus, the situation is quite similar to
that considered by Markov in~\cite{Markov1946}. The role of
the group $G$  from \cite{Markov1946}
is played by $G^*$, and the role of $A$ is played by $A^*$.
In \cite[Sections~3--5]{Markov1946}, Markov
defined sequences $\{b_i\}_{i=1}^\infty$,
$\{a_i\}_{i=1}^\infty$,  and $\{x_i\}_{i=1}^\infty$ with
certain properties (the properties of the first
sequence are labeled by 2.31 and 2.32 in his paper, the properties
of the second are labeled by 2.41--2.44, and the properties of the third
are labeled by 2.51 and 2.52).
We set $b_i= 1$ for all $i$; the sequences
$\{a_i\}_{i=1}^\infty$  and $\{x_i\}_{i=1}^\infty$ are already
constructed. These sequences have all
of the properties mentioned above except 2.32. However, 2.32 is not
used in Sections~6--12 of~\cite{Markov1946}; therefore,
the argument from these sections applies to the case under consideration
without any changes (except that $G$
should be replaced by $G^*$ and $A$ by $A^*$); in particular, the functions
$f_j$ on the sets
$$
A_j =\{1, a_j\}\cup \{a_k^{-1}x_ia_k:
i\ge j, k\le i\}
$$
given by
$$
f_j(1) = 0, \qquad f_j(a_j)= 1, \qquad f_j(a_k^{-1}x_ia_k) = \frac
1i
$$
are well defined and can be extended to seminorms $N_j$ on
the group $G^*$.  Following Markov~\cite[Section~11]{Markov1946}, for
any positive integer $n$ and sets of integers
$\{p_j\}_{j=1}^n$ and $\{q_j\}_{j=1}^n$, where $p_j\ge 0$ and $q_j>0$ for
$j = 1, \dots, n$, we define a seminorm $N^{p_1, \dots, p_n}_{q_1,
\dots, q_n}$ on $G^*$ by
$$
N^{p_1, \dots, p_n}_{q_1, \dots, q_n}(x) = \sum_{j=1}^n
N_{q_j}(a_{p_j}^{-1}x a_{p_j});
$$
we have $N_q^0 = N_q$ for $q= 1, 2, \dots$, because $a_0= 1$.
Markov showed that these seminorms determine a group
topology $\mathcal T^*$ on $G^*$; its neighborhood base
at the identity consists of all sets of the form
$$
U_N=\{x\in G^*: N(x) < 1\},
$$
where $N=N^{p_1, \dots, p_n}_{q_1, \dots, q_n}$ for some $p_1,
\dots , p_n\ge 0$ and $q_1, \dots, q_n>0$.

Let us prove that $1$ belongs to the closure of $A^*$ in this
topology. It is sufficient to show that $U_N$ intersects
$A^*$ for any $N$.

Let $N=N^{p_1, \dots, p_n}_{q_1, \dots, q_n}$. Take
a positive integer $s$ such that $s>n$, $s> p_j$, and $s> q_j$ for $j=1,
\dots, n$. By definition, we have
$$
N(x_s)=\sum_{j=1}^n N_{q_j}(a_{p_j}^{-1}x_s a_{p_j})= \frac ns < 1.
$$
Therefore, $x_s\in U_N$. On the other hand, $x_s\in A^*$ by construction.

Thus, the set $A^*$ is not closed in the group $G^*$ with the topology
$\mathcal T^*$. Since $G^*$ is a supernormal subgroup of $H$,
the neighborhoods of the identity in the topology $\mathcal T^*$
form a neighborhood base of the identity for some group
topology $\mathcal T$ on $H$. Since any neighborhood of the identity
in $\mathcal T^*$ intersects $A^*$ and $A^*=A\cap G^*\subset A$,
it follows that any
neighborhood of the identity in $\mathcal T$ intersects $A$. Thus, $A$
is not unconditionally closed in $H$.

Now, let $B$ be an arbitrary nonalgebraic set in $H$;
then $\widetilde B\ne B$. Take $b\in \widetilde B\setminus
B$. We have $1\notin b^{-1} B$. On the other hand, $1\in
b^{-1}\widetilde B$. According to
Lemma~12 from \cite{Markov1946}, $b^{-1} {\widetilde B}=
\widetilde{b^{-1}B}$; hence $1\in \widetilde{b^{-1}B}\setminus
b^{-1}B$. As was shown above, this implies the nonclosedness
of the set $b^{-1}B$ (and, therefore, of $B$ itself) in
some group topology on $H$. Thus, if the set
$B$ is not algebraic in $H$, then it cannot be
unconditionally closed in $H$. On the other hand, any algebraic
set in the group $H$ is unconditionally closed in
$H$~\cite[Theorem~1]{Markov1946}.
\end{proof}

\begin{corollary}
If $G$ is an Abelian group and $A\subset G$, then $A$ is
unconditionally closed in $G$ if and only if it is algebraic in $G$.
\end{corollary}

This assertion follows immediately from the fact that
any Abelian group can be embedded in a direct
product of countable groups as a subgroup (see, e.g.,
\cite{Kaplansky}) and that any subgroup of an Abelian group is supernormal
in this group.

\begin{corollary}
If $G$ is a direct product of countable groups and $A\subset G$, then
$A$ is unconditionally closed in $G$ if and only if it is algebraic in $G$.
\end{corollary}

Indeed, any countable subproduct in a direct product of groups
is a supernormal subgroup of the entire product.

Note that if $G_1$ and $G_2$ are two groups and $H_1\subset G_1$ and
$H_2\subset G_2$ are their supernormal subgroups, then, obviously,
$H_1\times H_2$ is supernormal in $G_1\times G_2$. This implies the following
assertion.

\begin{corollary}
If $G=G_1\times G_2$, where $G_1$ in an Abelian group and
$G_2$ is a direct product of countable groups, and $A\subset G$, then
$A$ is unconditionally closed in $G$ if and only if it is algebraic in $G$.
\end{corollary}

\end{document}